# Is incoherence required for sustainability?

Olivier Hamant[1,*]

[1] Michel Serres Institute and Laboratoire de Reproduction et Développement des Plantes, Université de Lyon, ENS de Lyon, UCBL, INRAE, CNRS, 46 Allée d'Italie, 69007 - Lyon, France

*Corresponding author
Correspondence: olivier.hamant@ens-lyon.fr

**ABSTRACT**

Unstoppable feedback loops and tipping points in socio-ecological systems are the main threats to sustainability. These behaviors have been extensively studied, notably to predict, and arguably deviate, dead-end trajectories. Behind the apparent complexity of such interaction networks, systems analysts have identified a small group of repeated patterns in all systems, called archetypes. For instance, the archetype of escalation is made of two positive feedback loops fueling one another and is prevalent when competition arises, as in arms race for instance. Interestingly, none of the known archetypes provide sustainability: they all trigger endless amplification. In parallel, in systems biology, there has been considerable advances on incoherent loops in molecular networks in the past 20 years. Such patterns in biological networks produce stability and a form on intrinsic autonomy for all functions, from circadian rhythm to immunity. Incoherence is the fuel of homeostasis of living systems. Here, I bridge both conclusions and propose that incoherence should be considered as a new operational archetype buffering socio-ecological fluctuations. This proposition is supported by the well-known trade-off between robustness and efficiency: adaptability requires some degree of incoherence. This applies to both technical and social systems: incoherent strategies recognize and fuel the diversity of solutions; they are the essential, yet often ignored, components of cooperation. Building on these theoretical considerations and real-life examples, incoherence might offer a counterintuitive, but transformative, way out of the great acceleration, and possibly, an actionable lever for decision makers.

*Keywords:* Incoherent loops, systems archetypes, systems thinking, robustness, sustainability, efficiency - resilience trade-off

# Introduction: Systems thinking with the help of archetypes

Human behaviors are often the result of untamed amplification loops. A bit like in a crowd movement: everyone is responsible for the agglomeration, and everyone wants to get out. Unfortunately, beyond a certain threshold, the crowd movement acquires its own logic, becomes autonomous, and is strengthened by the will of everyone to escape the crowd, paradoxically. Now considering global changes, one could say that humanity is currently in a planetary crowd movement: our socio-economic system is channeling our future trajectory through many amplifying loops. Alas, many of the proposed solutions do not fundamentally change this path. For instance, strategies relating to energy efficiency increments or to frugal smart technologies often ignore rebound effects, paradoxically leading to increased global consumption of resources and more global pollution in the end (1). Beyond energy, this applies to most of our products and activities, from concentrated detergents generating over-dosage to so-called smart technologies generating new needs in all sectors (2). How could one discriminate between these different options and identify truly sustainable strategies? Here I propose that incoherence is, paradoxically, the fuel of sustainability.

To develop this argument, one must dig into systems dynamics. As mentioned in *The limits to growth*, « running the same system harder or faster will not change the pattern as long as the structure is not revised » (3). This means that we need to deconstruct the interaction network behind the socio-ecological crisis and identify the amplification loops to be able to counteract them. This systemic work might seem difficult, but a guide exists, and is even very well known to systems analysts: it takes the form of "archetypes" (Figure 1).

An archetype is a small module, a sequence, that repeats itself in a system. If a sentence is composed of words, a system is composed of archetypes. More precisely, if a sentence has a meaning thanks to the syntax that links the words together, a system becomes predictable once the archetypes are articulated between one another. Such an approach is used widely, from the prediction of next day's weather to the multiple IPCC scenarios for the end of the century.

In the literature, there are about ten known archetypes only (4)(5), all describing amplification loops. For example, the "escalation" archetype best describes the crowd movement, as well as many other similar situations (Figure 1). Here is the case of intensive agriculture: ploughing, fertilizers and pesticides degrade soil structure and life; as agronomical yields decrease, more ploughing, fertilizers and pesticides are used to compensate for the decline in the short term, accelerating the agony of the soil in an endless escalation in the long term.

While the complexity of a system should not be reduced to its archetypes, the archetypes can be considered as the core skeleton of the system. This means that modifying the archetypes is, by far, the most effective way to change the trajectory of a given system. For instance, once the archetype behind traffic jam emergence is identified (the tragedy of the commons), one can find effective ways to prevent them, for instance, by imposing a slower speed on certain highway sections to enforce a more homogeneous behavior. Identifying the archetypes in a system is the first step to change its structure.

However, before starting this work, one must be sure that the list of known archetypes is complete. Here, I argue that we are missing one archetype, and arguably, the most relevant one to support sustainability: none of the canonic archetypes are stabilizing ones, whereas, in biological networks, there is increasing evidence that incoherent loops take on such a role. Should incoherence be considered as a new archetype in systems science, beyond biology? Could it become a relevant lever for sustainability?

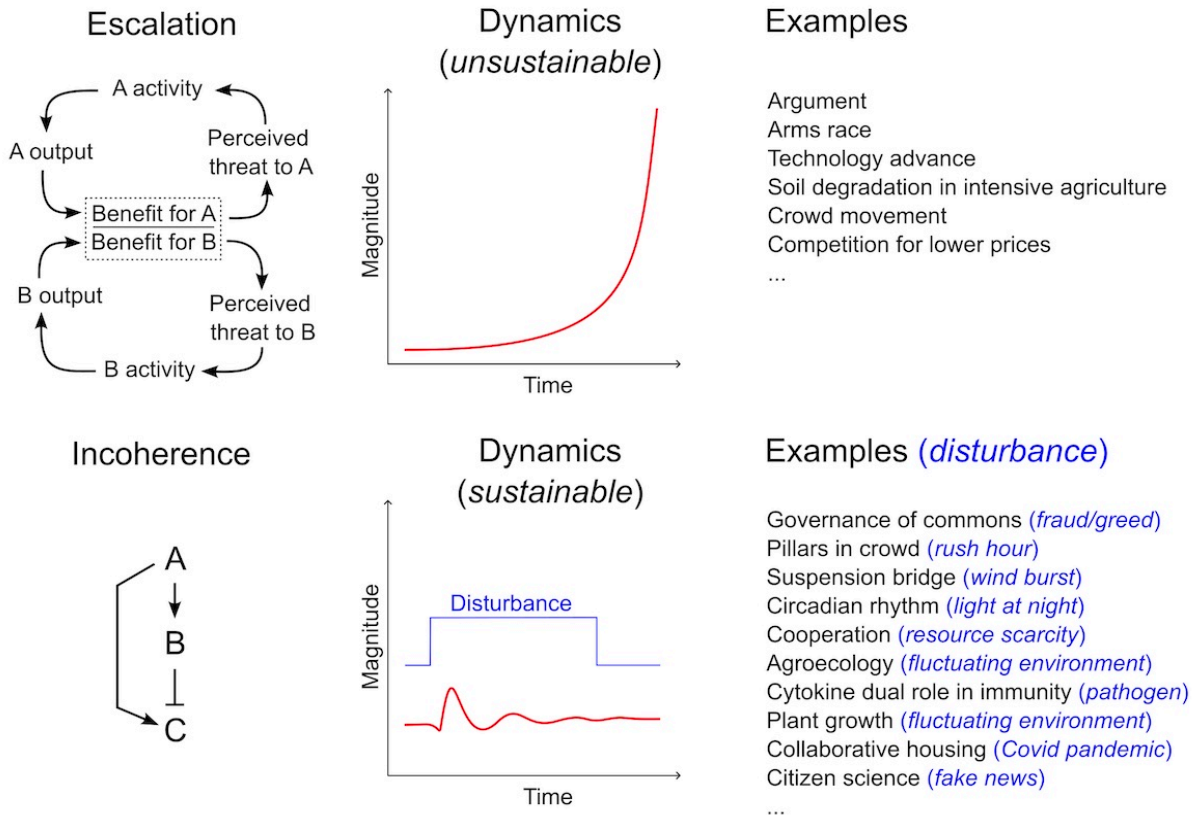

**Figure 1.** Incoherence as a sustainability systems archetype
Top: Escalation is a well-established systems archetype where A's output is perceived as threat by B, which in turn over-reacts, leading to a perceived threat from A, in a feedback loop(4)(5). The resulting dynamics is unsustainable, unless A or B drops the fight. Examples are listed on the right. Bottom: The example of the incoherent feedforward loop (well-established in biology). Such incoherence produces a form of robustness, the system becoming more resistant to internal or external disturbance. Many examples exist in biology (e.g. circadian rhythm (8)(9), immunity (20), or plant growth (21)). Incoherence is also an overarching factor in the governance of commons (16), addressing the shortcoming of the tragedy of the commons. Agroecology, citizen science or collaborative housing all exhibit incoherent features (slowness in a time of urgency, redundancy in a time of scarcity, heterogeneity in a time of rationalization) than more optimized and centralized approaches (intensive agriculture, academic/applied science, real-estate housing) but they are more robust to perturbations including environmental or economic fluctuations (through farmer's technical autonomy and biodiversity promotion in agroecology), fake news (through transdisciplinary debunking in citizen science) or psychological impact of Covid-imposed confinement (through increased interactions in collaborative housing). Incoherence's apparent weaknesses appear as systemic remedies to unsustainable amplifying loops.

**Incoherence is an essential ingredient of biological robustness**

Let us take some insights from systems biology first. That field has made significant progress in the past twenty years (6) and provides key conclusions: (i) amplification loops are over-represented in biological networks (gene networks, but also neural networks, ecosystem networks, etc., gene networks being the most exhaustively studied), (ii) these loops are almost always built in such a way that one branch has the opposite effect of the other, i.e. they are incoherent (Figure 1). Interestingly, these incoherent loops have three main effects: (i) they allow the system to return to its original stable state after a perturbation, (ii) they are rather inhibitory to avoid endless exponential amplifications, (iii) they improve very significantly the inhibition dynamics, and thus the induced stability (7).

More simply said, incoherent loops produce stable outcome, through autonomous oscillatory behavior that also dampen external perturbations. This is for example the case of circadian rhythm in living beings: (i) an amplification loop in the molecular network generates two possible outputs (active metabolism during the day and slow metabolism at night); (ii) an inhibition loop generates oscillations; (iii) when combined, they generate an oscillating bistability (day/night metabolism in alternation). This behavior is very robust in the long term, i.e. little disturbed by external elements, thanks to the incoherent loop structure of the system. This type of behavior is now formally established in mammals, insects, fungi or plants (8)(9).

Many other examples can be found in biological systems. For instance, plant growth is driven by their internal hydrostatic pressure. This is obvious when one stops watering them: plants wilt, i.e. they deflate. When they are inflated, they grow. However, in these conditions, plants are using photosynthesis to build stiff walls around their cells, to resist their own turgor pressure. In other words, plants drive their growth with the handbrake on (10). This incoherence is a key to their robustness: many studies show that this constant balance allow them to constantly adapt to their environment, and to better use external and internal fluctuations to shape their organs (11).

Such biology insight is an opportunity to revisit the definition of stability though systems dynamics: stability is not a flat line; it is instead a dynamic oscillation, alternating action and reaction. In other words, the autonomous oscillation of the system can protect against unforeseen and random fluctuations. For instance, the cells of the gut are renewed periodically (in a few days) and this allows resistance to a great number of external aggressions. Conversely, the non-renewal of teeth limits the range of possibilities. This also echoes Ashby's law of requisite varieties, where systems robustness does not emerge from its further consolidation (i.e. increasing coherence), but instead from its own variability (i.e. fueling incoherence) (12). While coherence feeds amplifying loops and generate instability, incoherence-fueled oscillation is an effective shield against turbulence in a fluctuating environment.

## Incoherence fuels robustness in technical systems too

Since we are dealing with systems, and therefore only with logical links, this conclusion also applies to other contexts. For instance, there is strong interest in engineering for these incoherent loops, in particular to avoid uncontrollable exponential behavior. Beyond elaborate examples in electronics or robotics, one can think of the simpler case of suspension bridge. It displays a fundamental contradiction: the bridge piers are in compression while the deck cables are in tension. Two destructive forces in a bridge, for what benefit? Because they oppose each other, the forces - compression and tension - create a mechanical balance. This allows the bridge to have its own mechanical autonomy, and thus to resist external fluctuations: when subjected to wind, the bridge oscillates but does not break, at least up to a certain deformation threshold. Needless to say, such a balance of force also exists in (robust) biological systems, and typically in cells, with a membrane under tension and a content under compression (13).

Could incoherence, as a stabilizing archetype, help us design sustainable strategies counteracting the amplification loops of the Anthropocene? Let us consider the rise of low tech. While the world is experiencing more pressing issues, this could be viewed as an incoherent idea: creating suboptimal solutions, through slow participatory research with citizens, and with reduced economic growth (e.g. the all-repairable), instead of favoring the most optimized solution. The contradiction is resolved when considering that the race toward optimization also makes our socio-economic and ecological systems increasingly vulnerable. Robustness emerges from suboptimal processes, and prevails over efficiency, in a fluctuating environment (14). Optimization is relevant to a stable environment, only.

Could incoherence feed sustainability beyond technical solutions? In an increasingly fluctuating world, the key, shared, values is no longer efficacy or efficiency, it is adaptability (14). Can incoherence be used as a tool to stabilize our turbulent world?

## Incoherence fuels robustness in organizations too

Incoherence, as an effective lever for sustainability, also applies to organizations. Consider the archetype of "the tragedy of the commons": each individual pursues their own goals in a selfish manner, and in the end, the commons (e.g. in pastoral areas: a shared aquifer, a shared grazing plot or a shared irrigation system (Figure 2) no longer fulfill their function to the detriment of all (15). The principles that prevent the emergence of the tragedy of the commons have been identified by economist Elinor Ostrom and her team, through a global analysis of robust commons over time (16) (Figure 1). Most of these principles are fundamentally incoherent: (i) a commons which primary value could be its openness, instead needs limited access; (ii) the rules to maintain a fragile resource could be strict, but the analysis of robust commons over time shows instead that these rules must be modifiable in a participatory manner; (iii) sanctions, which one might want very strong to limit selfish behavior and fraud, must on the contrary be weak to guarantee the cohesion and belonging of group members; (iv) commons survive if they are managed by a self-organized and autonomous structure, but this structure must be recognized by an external entity in order to be maintained over time. Thus, it appears that incoherence might be a necessary, sometimes counterintuitive, ingredient to correct amplification loops, even in organizations.

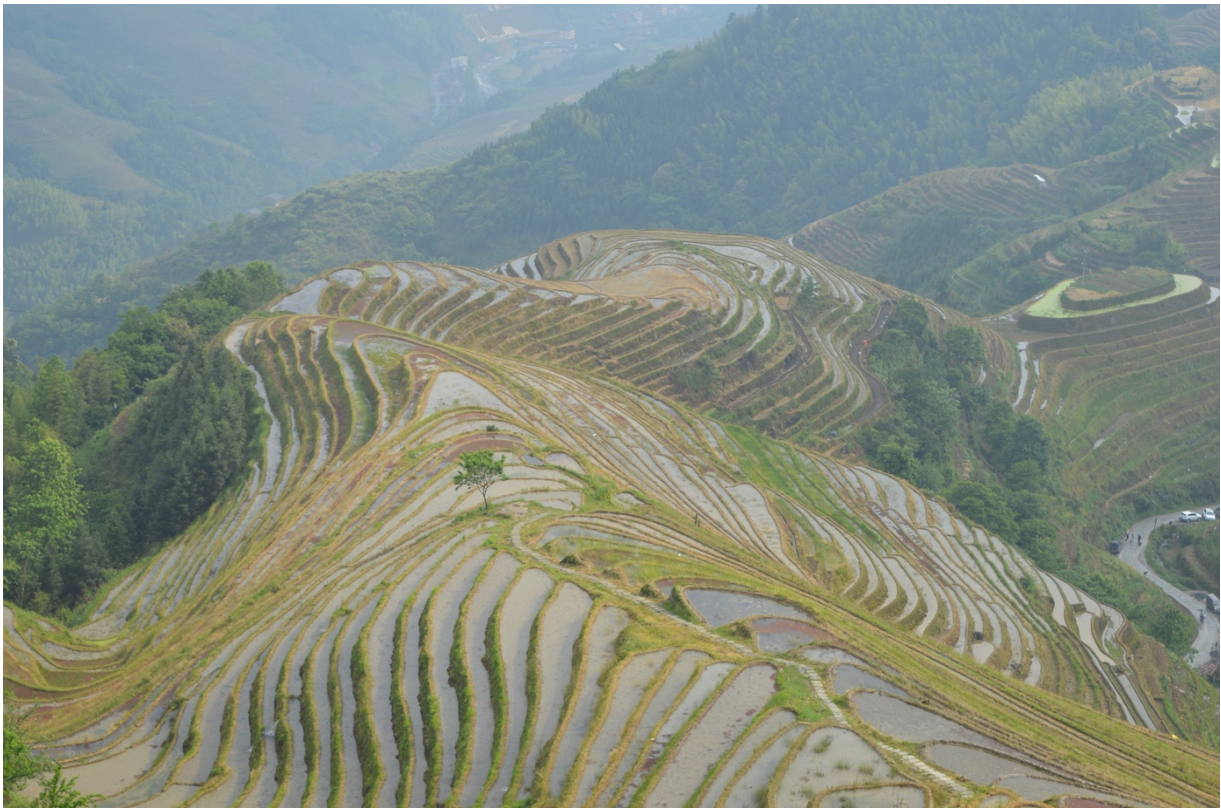

**Figure 2.** A rice field in China, as an example of sustainable system requiring the governance of the commons, and associated incoherence.
Coordination between upstream and downstream farmers requires cooperation, i.e. a shared objective that may conflict short-term immediate interests. Elinor Ostrom's principles from the governance of the commons embed several incoherent components, which in the end, ensure the sustainability of the resource. See main text for details. Credit photo: O. Hamant.

To take another example, the nurse collective Buurtzorg (in the Netherlands) has identified a solution to the burnout and "big quit" epidemic in the health sector through a counterintuitive solution: they gave up commercial objectives. Spending more time with patients initially triggered increased expenditure, but, in the end, this made patients more autonomous and thus better taken care of, and nurses more aligned with their work values and more engaged. This also led to positive financial return to the State, ultimately. Now Buurtzorg gathers ca. 30'000 nurses, also demonstrating how this counterintuitive model can scale up (17). Such a strategy echoes the work on loosely coupled systems, which are well-known for their incoherence and robustness (18). Robustness does not emerge from tighter interactions, but instead from increased noise and more oil in the wheels.

Another incoherent strategy is to degrow and disperse to survive (instead of growing further to become stronger). This echoes the "r-strategy" in living beings, and notoriously in flowering plants: making many seeds, and dispersing them to ensure the survival of the species against all odds. In social systems, this can take the form of decentralization, e.g. selling local milk bricks, at a fair price in consultation with farmers, and outside of large central purchasing agencies. A decentralized link to territories is less efficient and may seem wasteful from the lean management perspective, but it has a greater chance of respecting citizens and ecosystems, precisely because it holds its own limitation in its structure.

## Conclusion: incoherence is required for cooperation

More generally, the power of incoherence can be illustrated with the difference between collaboration and cooperation. When people collaborate, they each move forward on their respective projects, hoping that the addition of these contributions will be synergistic. However, a successful collaboration often masks destructive effects for unidentified third parties. Collaboration can thus hide a form of competition. On the contrary, when people cooperate, the success of the group project outweighs the success of individual projects. Cooperation therefore implies an acceptance of inconsistency, even to the point of having to score an own goal. Choosing cooperation over collaboration means understanding, as in Michel Serres' Natural Contract (19), that any transaction between humans involves multiple social and environmental interactions, including incoherent ones, that go beyond short-term interests, for the survival of the group.

With these examples, I hope to raise awareness in sustainability science and policy: incoherence, bridging recent advances in biology network analysis and established framework in systems science, can provide a way out of amplifying loops. Beyond the examples mentioned above, one can evoke the amplifying loops generated by social networks and artificial intelligence, through which the same information is distributed based on the coherent profiles of users. To break such regimentation, one needs to depart from the shared beliefs in the group, and thus own some degree of incoherence.

In the end, incoherence might be the most relevant nature-inspired innovation for sustainability. Going back to the opening example, to prevent crowd movement, one can install fixed pillars on the walking path beforehand to slow down everyone at rush hour, another incoherent solution (allowing people to escape by slowing them down). Incoherence provides an answer to our over-optimized, and thus vulnerable, civilization, building on apparent inefficiencies (heterogeneity, slowness, errors, conflicts, randomness, etc.). It is an effective remedy to unstoppable amplifying loops and an actionable tool for decision-makers.

## Acknowledgements

I am thankful to Dennis Meadows, the Balaton group and the Michel Serres Institute for helpful discussions on this topic.

## Conflict of interest disclosure

The authors declare that they comply with the PCI rule of having no financial conflicts of interest in relation to the content of the article.

## Funding

The author declares that he has received no specific funding for this study.